\documentclass[12pt]{article}
\usepackage{epsfig}
\usepackage{amsfonts}
\def\R{{\mathbb R}}
\def\C{{\mathbb C}}
\def\L{{\mathcal L}}
\def\sbp{\subparagraph}
\title{Diagrams of divide links}
\author{Sergei Chmutov}
\date{}

\begin{document}


\maketitle

 There are several
fascinating relations of plane immersed curves and links. One of
them which goes through Legendrian links led Arnold \cite{Ar} to
discovery of three simple invariants $J^+$, $J^-$, $St$ of such a
curve. N.~A'Campo in \cite{AC-1} suggested another construction of
a link from a generic immersion of a curve into a 2-disk. It is
tightly related to the singularity theory. If $f_\R: (\R^2,0)\to
(\R,0)$ is a germ of an analytic function whose complexification
$f_\C$ has an isolated critical point then one can define a link
$\L$ of the singularity $f_\C$ as an intersection of the
zero-level variety of $f_\C$ with a small 3-sphere in $\C^2$
centered at the critical point. Such links are sometimes called
{\it algebraic links}. On the other hand, in singularity theory it
is useful \cite{AC-0, GZ} (see also \cite{AGV}) to consider a
small perturbation $D$ of the real plane singular curve
$\{f_\R=0\}$ which is a generic immersed curve with the maximal
possible number of double points. A'Campo \cite{AC-1} restored the
link $\L$ directly from the curve $D$.

In this paper we give a simple procedure to draw a diagram of the
link $\L$ from a picture of $D$ (Theorem 2.2). It is essentially a
particular case of the results of \cite{CP}. An advantage of our
approach is that we obtain a link diagram directly from a divide
picture without deforming it into so called ordered Morse signed
divide as in \cite{CP}. A similar method to draw diagrams was
found by M.~Hirasawa in \cite{Hi}. Our diagrams are obviously
symmetrical in a sense that the rotation of $\R^3$ by $180^\circ$
around the $x$-axis reverses the orientation of $\L$.

\section{Divides and their links.}

\sbp{1.1. Definition}
\label{de-1.1}
(\cite{AC-1, AC-2}). {\it A divide} $D$ is the image of a generic
immersion of a finite number of copies of the unit interval
$I$=[0,1] in the unit disk $B\subset\R^2$ such that $\partial I$
is embedded in $\partial B$ and double points are the only
singularities allowed.

We consider divides up to isotopy of the disk $B$. The isotopy does
not assume to be identical on the boundary $\partial B$.

\sbp{1.2. Example.}
\label{ex-1.2}
The curve $x^3+y^4=0$ has a singularity of type $E_6$ at the
origin \cite{AGV}. A small perturbation of it is a divide which
looks as follows.
$$
  \mbox{\begin{picture}(60,60)(0,0)
  \put(0,0){\epsfxsize=60pt \epsfbox{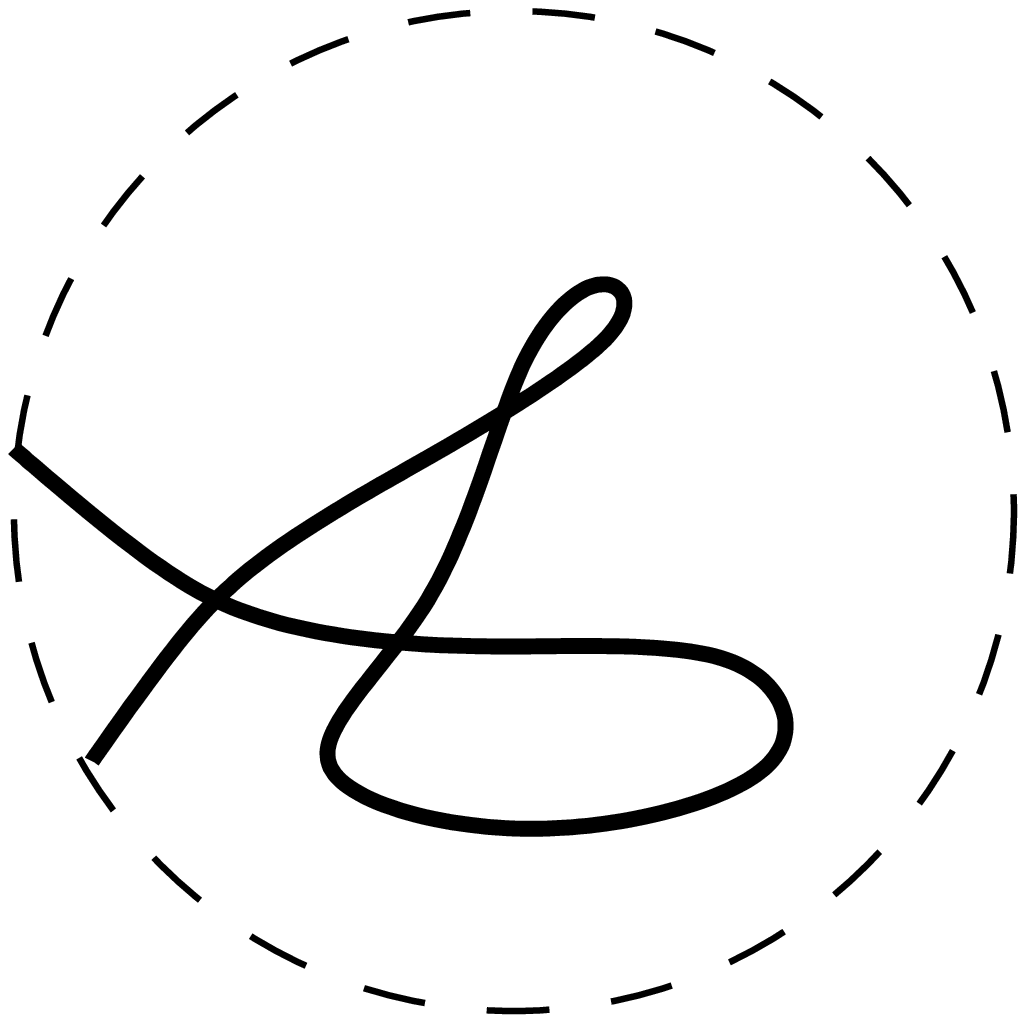}}
        \end{picture}}
$$

\sbp{1.3. Definition} (\cite{AC-1, AC-2}). \label{de-1.3}
Let $x$ be the horizontal coordinate on the disk $B$ and $y$ be
the vertical coordinate. {\it A divide link} $\L_D$ is a link in
the 3-sphere $S^3=\{(x,y,u,v)\in\R^4 | x^2+y^2+u^2+v^2=1\}$ such
that $(x,y)$ is a point on $D$ while $u$ and $v$ are the
coordinates of a tangent vector to $D$ at the point $(x,y)$.

So each interior point of $D$ has two corresponding points on $\L_D$,
and a boundary point of $D$ gives a single point on $\L_D$.

The link $\L_D$ has a natural orientation. Indeed, choose any orientation
of every branch of $D$. Let $(u,v)$ be the tangent vector to $D$ at
$(x,y)$ pointing to the direction of the chosen orientation of $D$. Then
the orientation of $\L_D$ is given by the vector
$(\dot{x},\dot{y},\dot{u},\dot{v})$.
It is easy to see that this orientation of
$\L_D$ does not depend on the choice of orientations of branches of $D$.

The number of components of $\L_D$ equals to the number of
branches of the divide $D$ which is the number of the copies of
the unit interval $I$ in the definition \ref{de-1.1}. In
particular, if $D$ consists of only one branch like in
\ref{ex-1.2} then $\L_D$ will be a knot.

\sbp{1.4.}
 Topological type of a divide link does not change under
a regular transversal isotopy of the disk $B$. So it does not
depend on the choice of coordinates in \ref{de-1.3}. Also it does
not change under a moving of a piece of the curve $D$ through a
triple point \cite{CP}. In particular, the following two divides
have the same knot type as one in 1.2.
$$
  \mbox{\begin{picture}(300,60)(0,0)
  \put(0,0){\epsfxsize=60pt \epsfbox{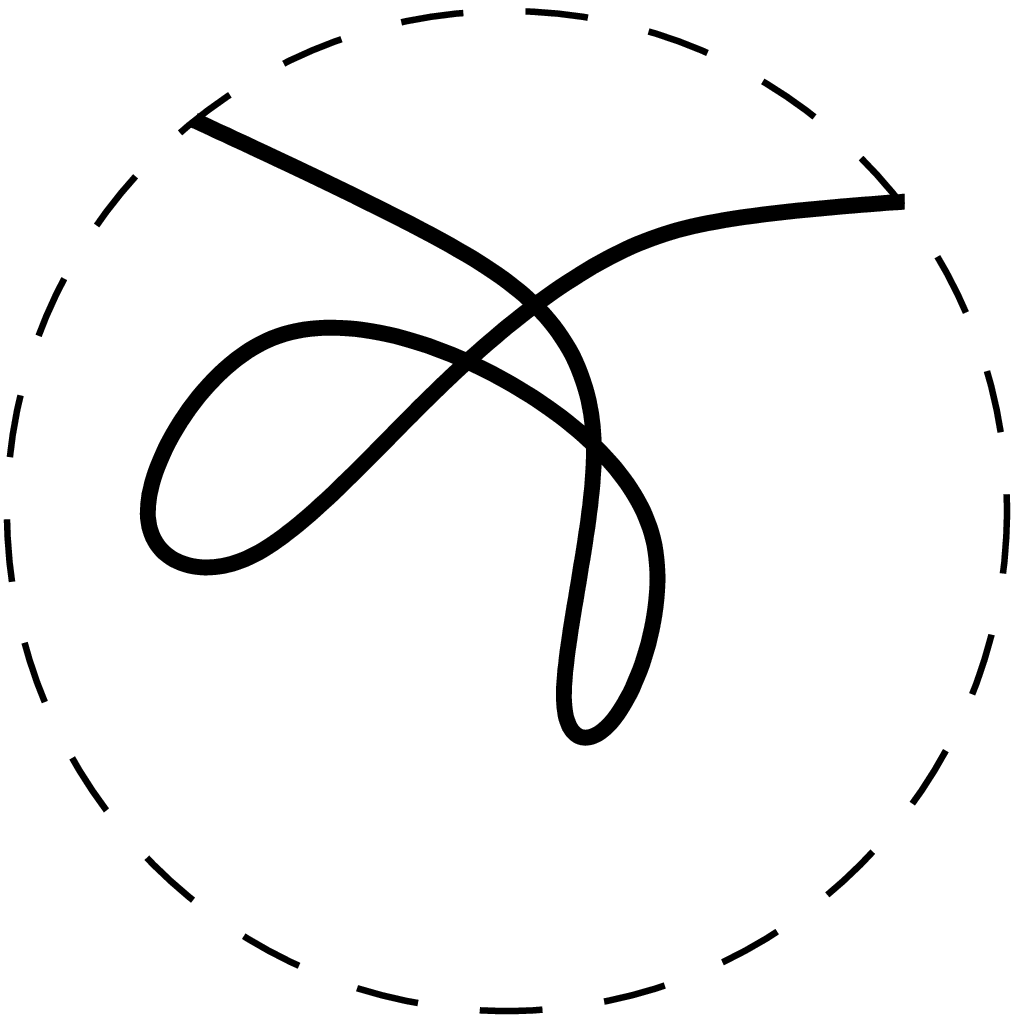}}
  \put(200,0){\epsfxsize=60pt \epsfbox{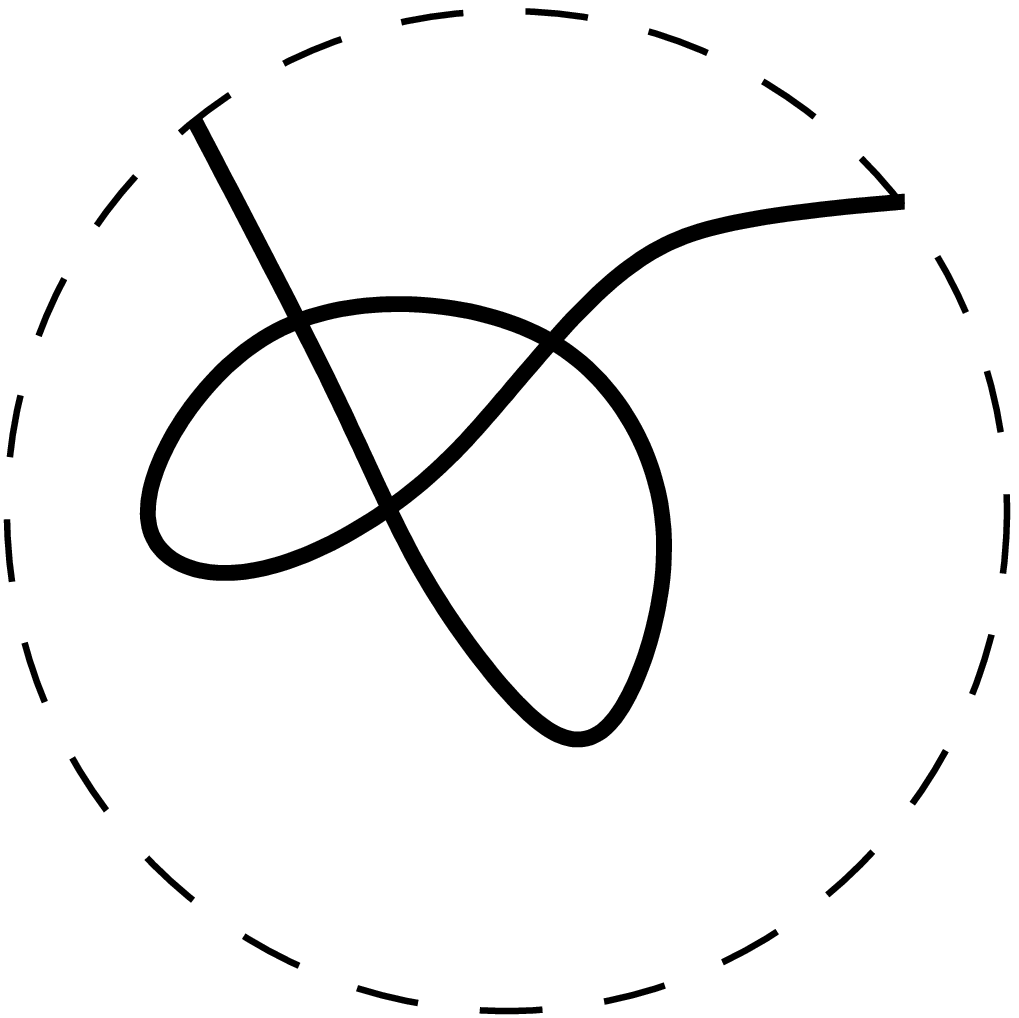}}
        \end{picture}}
$$

\sbp{1.5.}
In \cite{AC-1} A'Campo proved that all algebraic links 
are
divide links. In \cite{AC-2} he showed that the links $\L_D$
corresponded to a connected divide $D$ are fibered and computed
their monodromy in terms of the combinatorics of divide $D$.
Moreover, he proved that the unknotting number of a one-branch
divide knot $\L_D$ equals the number of double points of $D$. Not
all fibered links have the form $\L_D$. Figure eight knot $4_1$ is
not a divide knot. It is not clear how large is the class of
divide links in the class of all fibered links.

\sbp{1.6.}
It is known for a long time \cite{Bu} that algebraic 
knots are
classified by the Alexander polynomial. N.A'Campo (\cite{AC-3})
found two different divide knots with the same Alexander
polynomial.

\section{Diagrams of divide links.}

\sbp{2.1.}
Let us call a divide {\it generic} if its points with
vertical tangents differ from the double points and the boundary
points and their $x$-coordinates are pairwise different. The
divide in \ref{ex-1.2} is generic. Any divide can be made generic
by a small deformation.

\sbp{2.2. Theorem.} 
\label{th-2.2}
{\it
For a generic divide $D$ a link diagram of
$\L_D$ can be drawn in the following way:\vspace{-5pt}
\begin{enumerate}
\item Consider a horizontal line below the disk $B$, say the line
    $\{y=-1.5\}$. Let $s$ be the symmetry  with respect to the line. We
    are going to draw the link diagram of $\L_D$ by modifications of
    the union $D\cup s(D)$.
\item Replace each double point of the union by a crossing of the type
  \quad $\mbox{\begin{picture}(20,10)(0,0)
        \put(0,-2){\epsfxsize=20pt \epsfbox{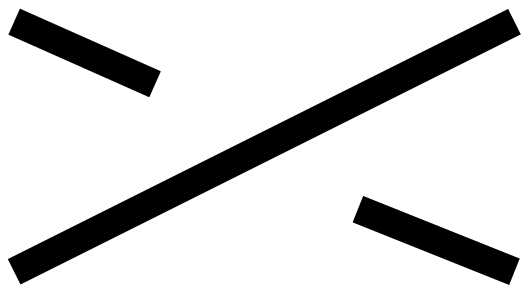}}
        \end{picture}}$\ .
\item Connect each boundary point $p$ of $D$ with $s(p)$ by a
vertical string.
\item Replace a small piece of our curve near each point $p$ with
vertical tangent by two vertical strings connecting the upper and
lower parts of the picture as in examples below. The strings make
a positive half-twist at the line of the reflection $s$. Note that
at the upper part of the picture the right string goes below all
intersected intervals of $D$ while the left string goes above the
intervals. Correspondingly at the lower part of the picture the
right string goes above the intervals of $s(D)$ and the left
string goes below the intervals. The examples \ref{ex-2.3} and
\ref{ex-2.4} demonstrate this.
\end{enumerate}
}

\sbp{2.3. Example.} \label{ex-2.3}
For the divide \ref{ex-1.2} the theorem gives the following
diagram:
$$
  \mbox{\begin{picture}(360,250)(0,0)
  \put(-5,165){\epsfxsize=60pt \epsfbox{e6-0.eps}}
  \put(71,195){\epsfxsize=25pt \epsfbox{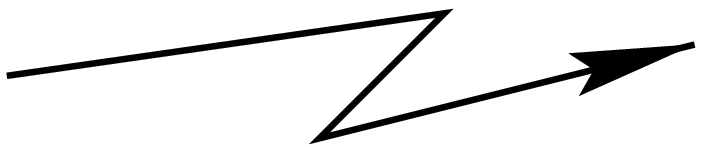}}
  \put(75,210){\mbox{(1)}}
  \put(110,150){\epsfxsize=80pt \epsfbox{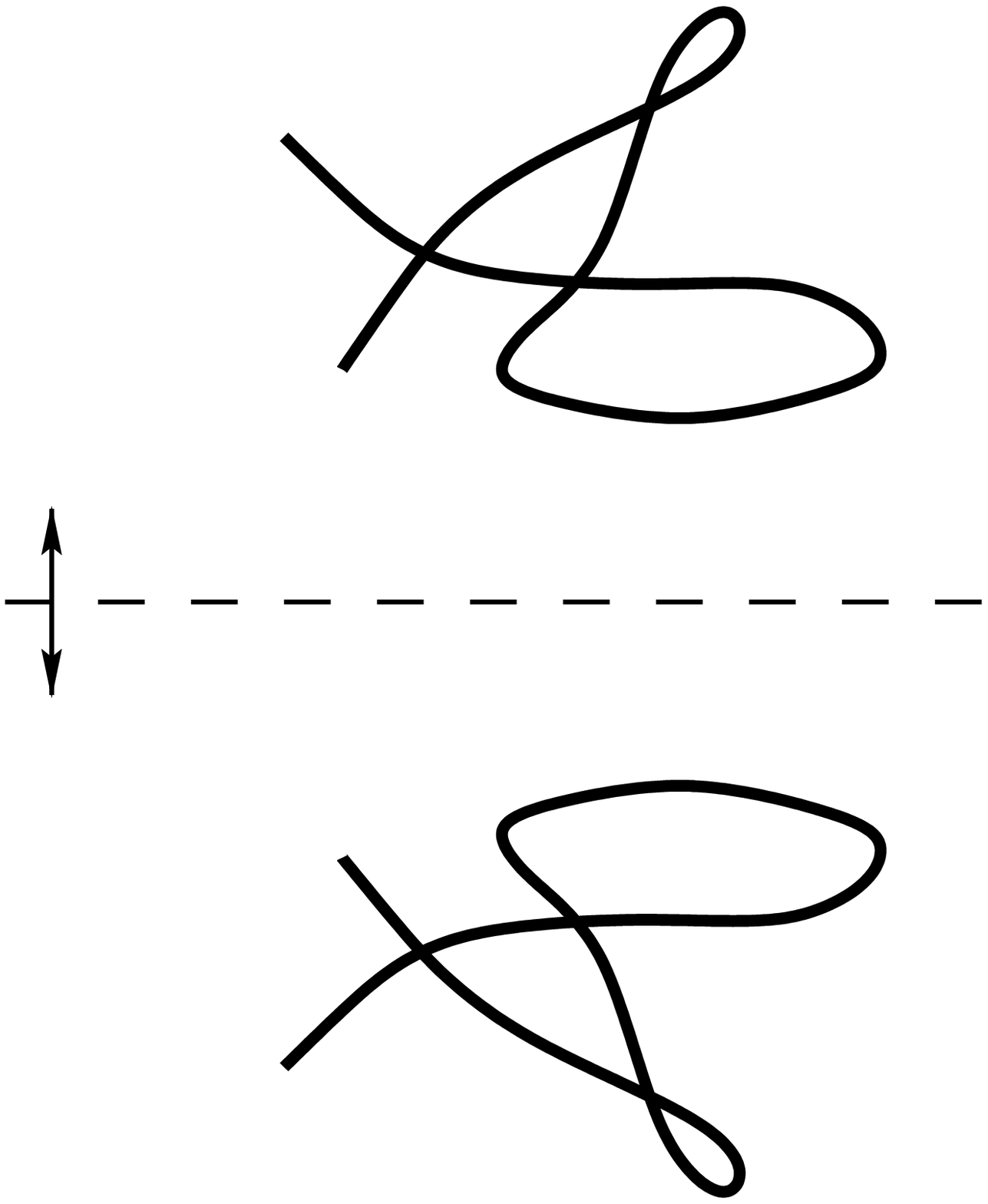}}
  \put(115,183){\mbox{$s$}}
  \put(135,240){\mbox{$D$}}
  \put(135,150){\mbox{$s(D)$}}
  \put(195,198){\mbox{$\{y=-1.5\}$}}
  \put(261,195){\epsfxsize=25pt \epsfbox{toto.eps}}
  \put(265,210){\mbox{(2)}}
  \put(300,140){\epsfxsize=60pt \epsfbox{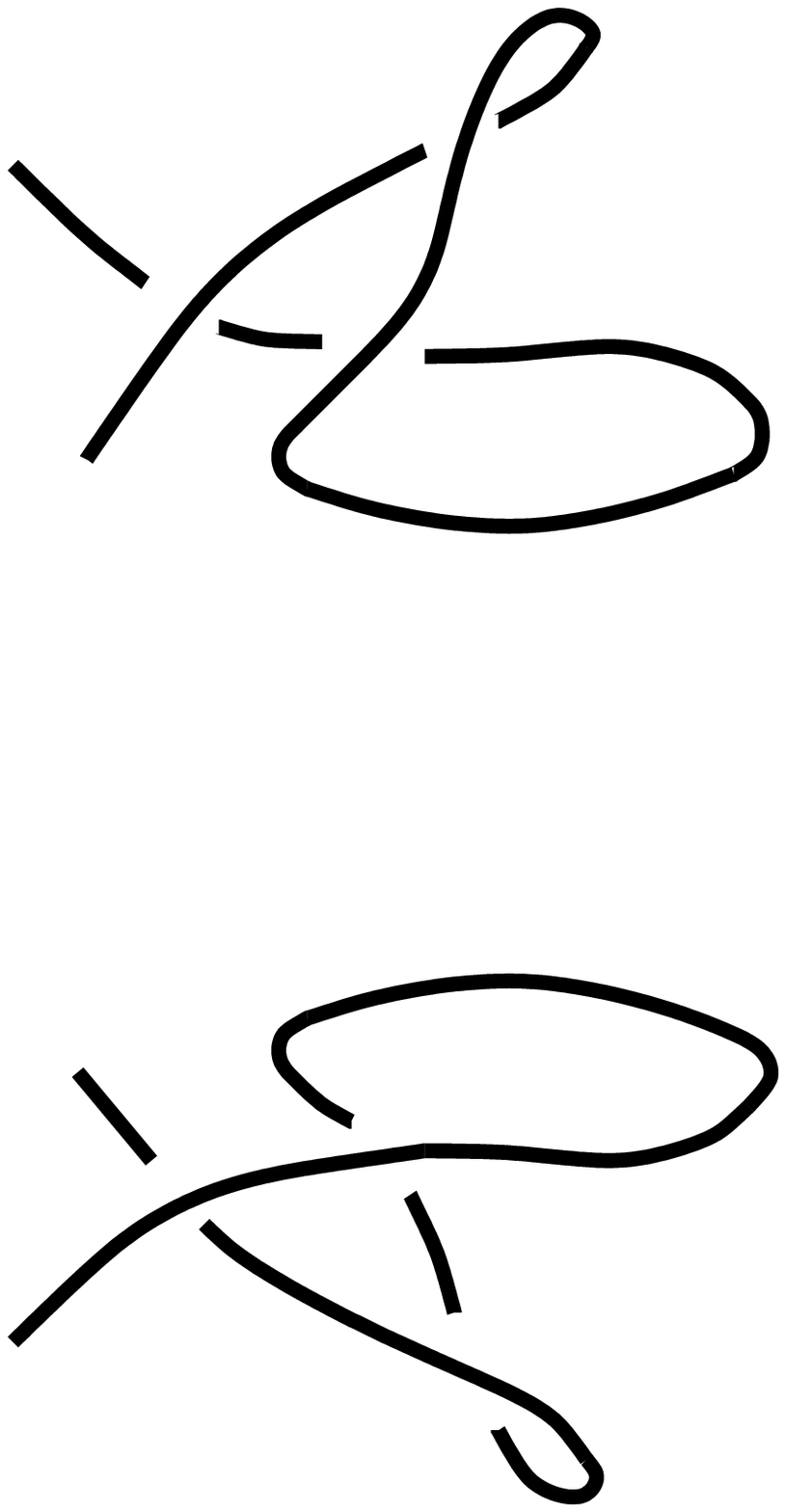}}
  \put(76,55){\epsfxsize=25pt \epsfbox{toto.eps}}
  \put(80,70){\mbox{(3)}}
  \put(120,0){\epsfxsize=70pt \epsfbox{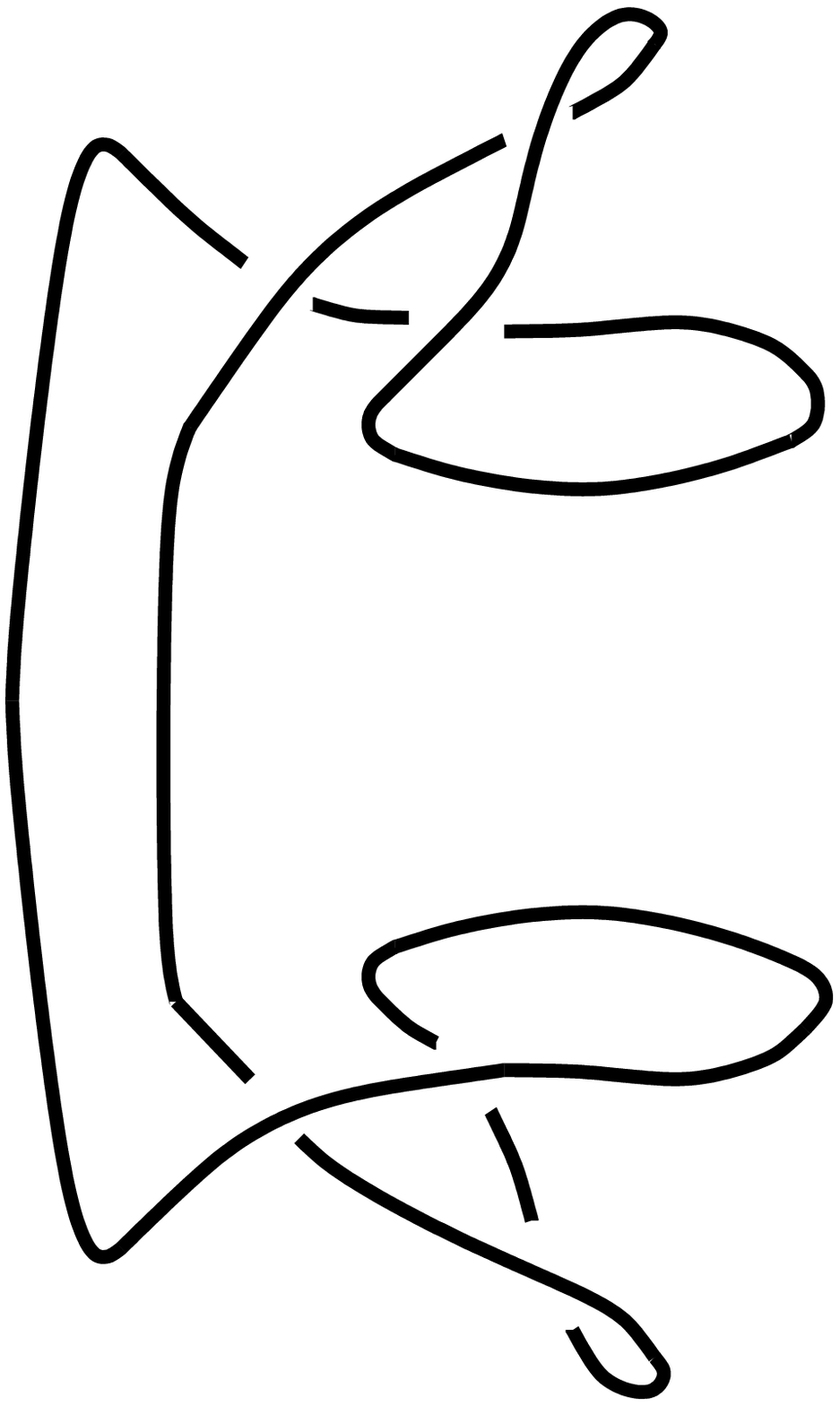}}
  \put(216,55){\epsfxsize=25pt \epsfbox{toto.eps}}
  \put(220,70){\mbox{(4)}}
  \put(260,0){\epsfxsize=75pt \epsfbox{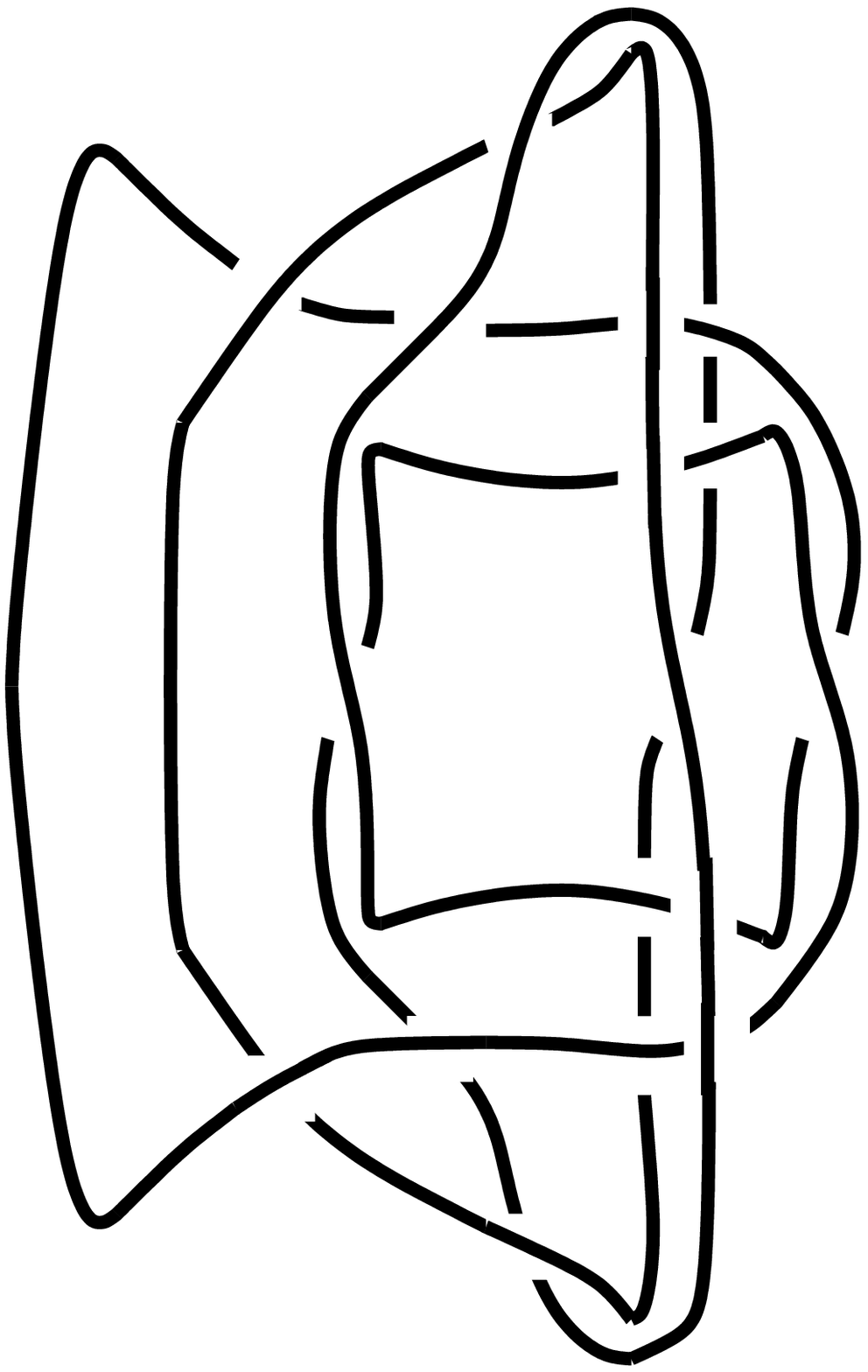}}
        \end{picture}}
$$

\sbp{2.4. Example.}\label{ex-2.4}
We borrow this divide from \cite{AC-2}. The corresponding knot
is 10$_{145}$.
$$
  \mbox{\begin{picture}(300,100)(0,0)
  \put(30,15){\epsfxsize=60pt \epsfbox{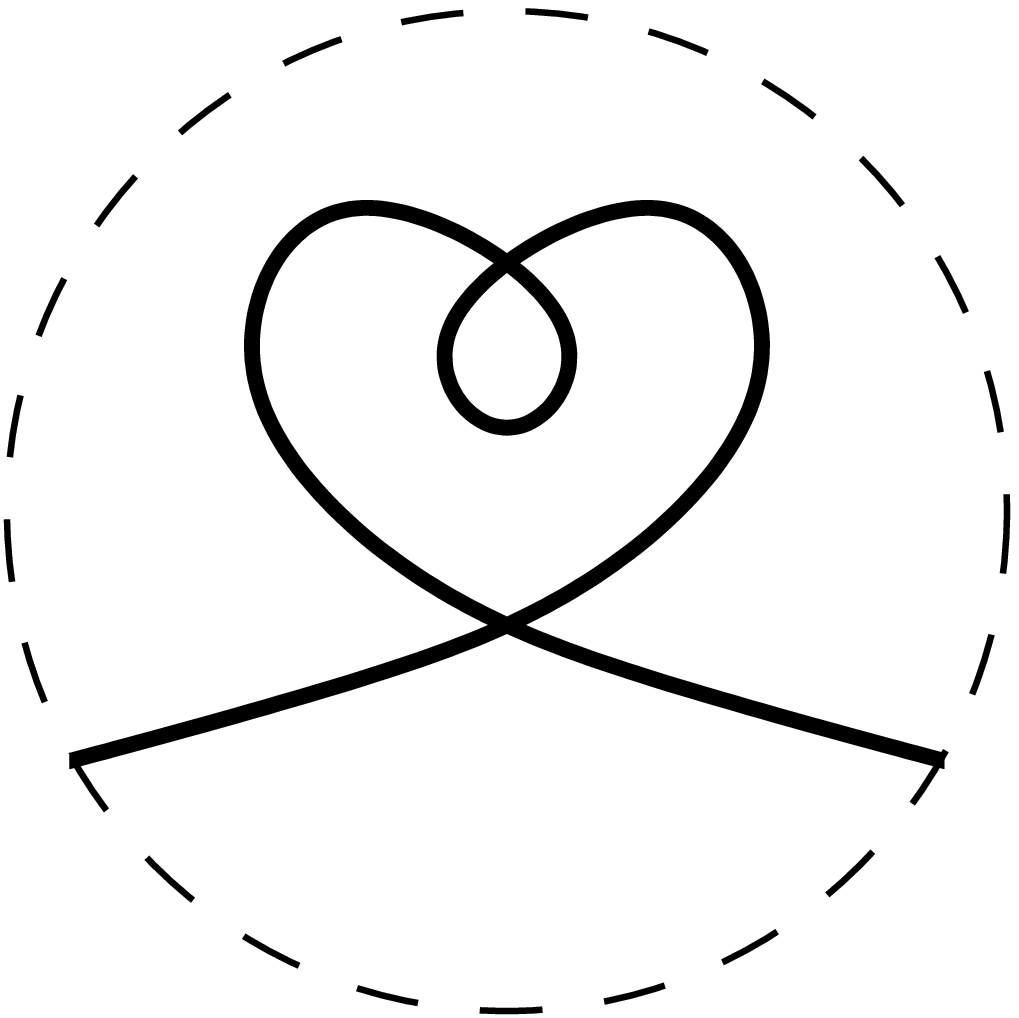}}
  \put(120,40){\epsfxsize=25pt \epsfbox{toto.eps}}
  \put(170,-15){\epsfxsize=75pt \epsfbox{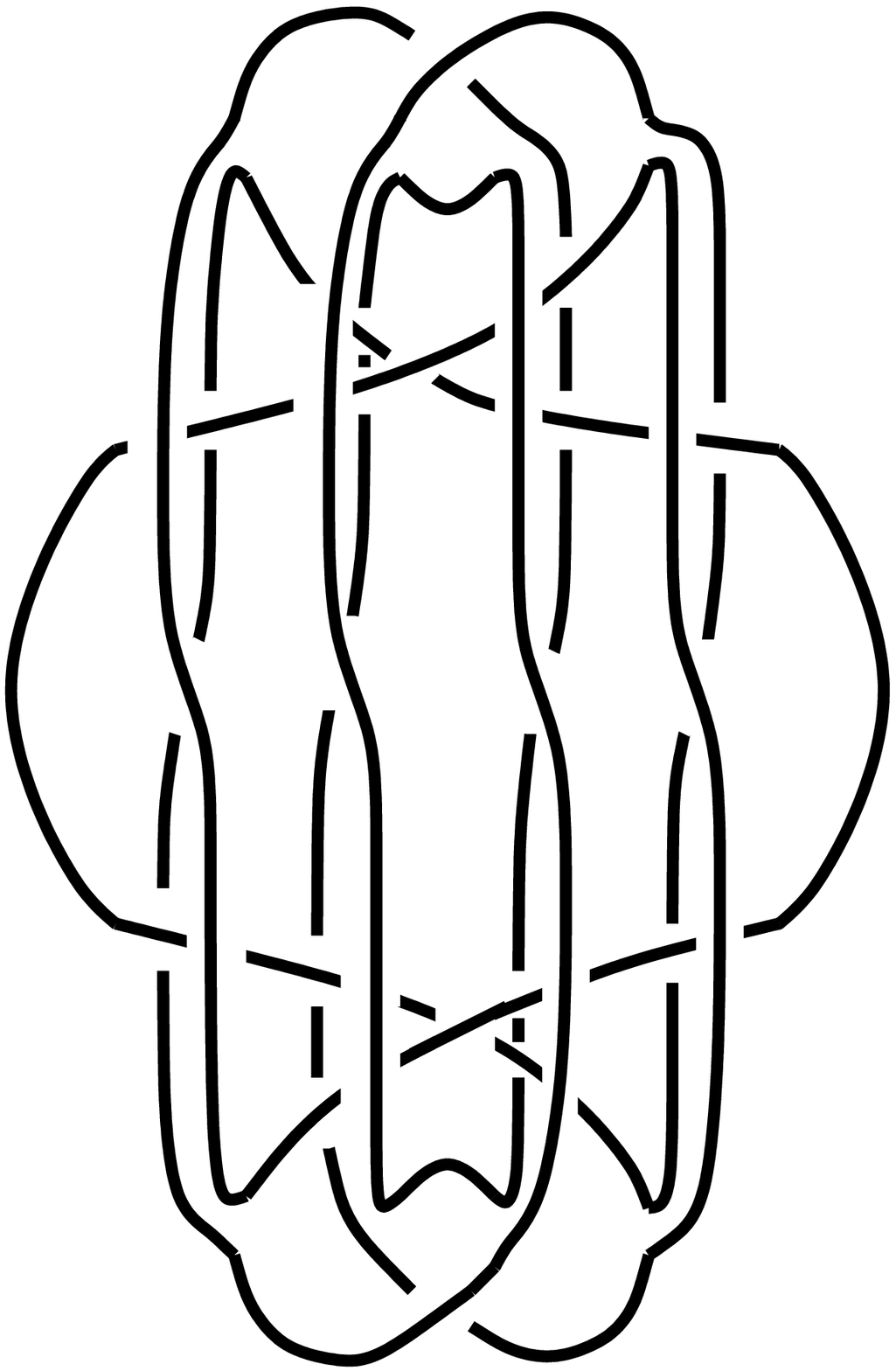}}
        \end{picture}}
$$

\bigskip\bigskip
\sbp{2.5. Remark.} 
As it was noted in \cite{CP} the theorem is also valid in the
situation when we allow closed immersed components in the
definition of a divide.

\sbp{Proof of the Theorem \ref{th-2.2}} The theorem follows from
the results \cite{CP} which give a representative braid for {\it
ordered Morse signed divide} (OMS). An OMS is a divide $D$ such
that $x$-coordinate as a function on $D$ has only two critical
values: $a$ as the minimum critical value and $b$ as the maximum
critical value; and $x$-coordinate of each double point is between
$a$ and $b$. Besides this a sign $+$ or $-$ is attached to each
double point of $D$. See the details in \cite{CP}. An OMS is not a
generic divide in our sense. Also in \cite{CP} there is an
algorithm transforming an arbitrary divide to OMS (and attaching
signs to double points). We use a particular such transformation
pulling down a narrow tail near each critical point of the
function $x|_{_D}$ and then move a minimum (resp. maximum) left
(resp. right) to the level $a$ (resp. $b$). For the divide of
example \ref{ex-2.4} this transformation gives the following OMS:
$$
  \mbox{\begin{picture}(350,150)(0,0)
  \put(50,40){\epsfxsize=60pt \epsfbox{ex24.eps}}
  \put(140,75){\epsfxsize=25pt \epsfbox{toto.eps}}
  \put(170,4){\epsfxsize=150pt \epsfbox{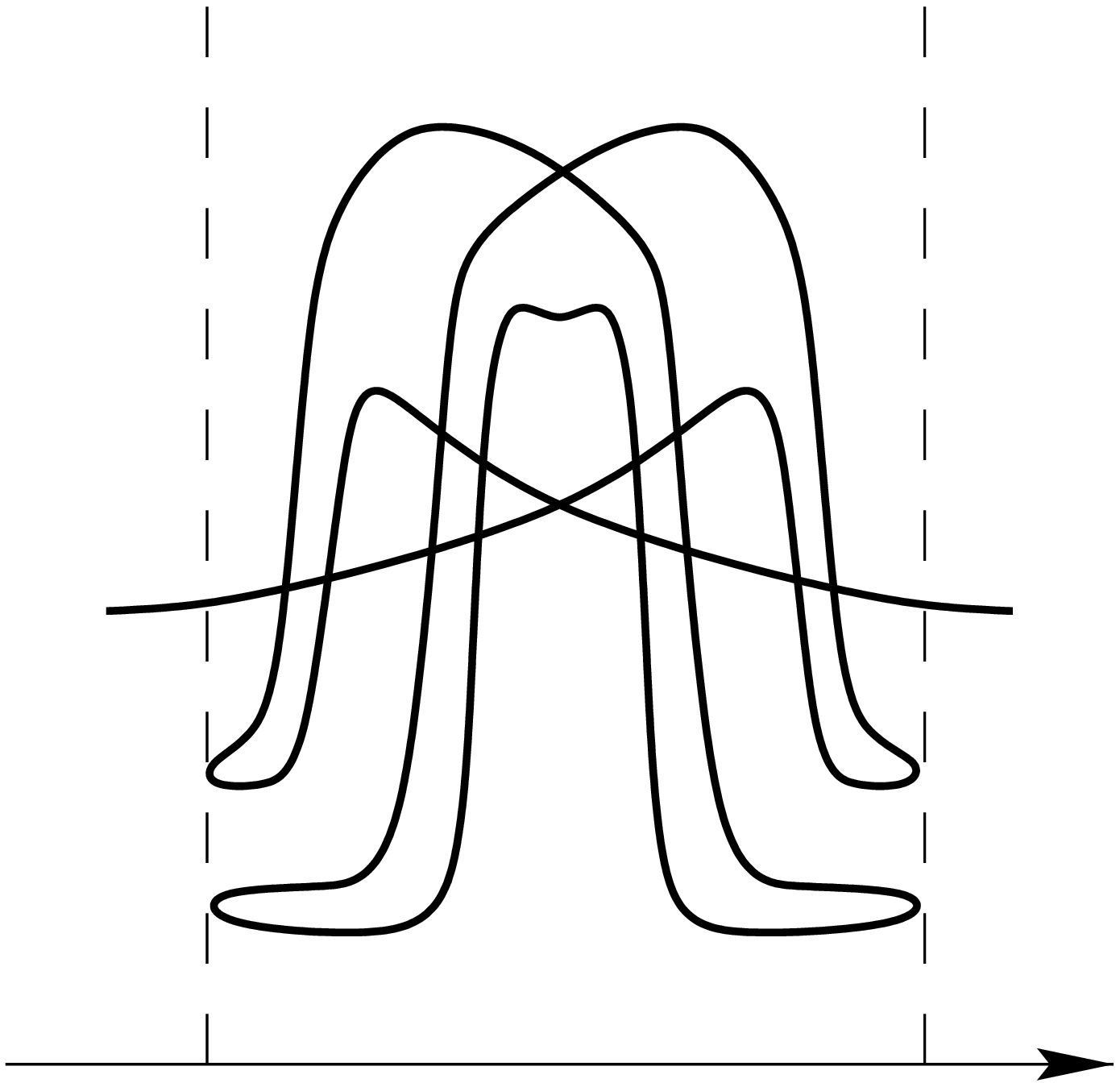}}
  \put(200,75){\mbox{$\scriptstyle +$}}
  \put(216,69){\mbox{$\scriptstyle -$}}
  \put(221,79){\mbox{$\scriptstyle +$}}
  \put(236,74){\mbox{$\scriptstyle -$}}
  \put(283,75){\mbox{$\scriptstyle +$}}
  \put(269,68){\mbox{$\scriptstyle -$}}
  \put(263,79){\mbox{$\scriptstyle +$}}
  \put(248,74){\mbox{$\scriptstyle -$}}
  \put(242,76){\mbox{$\scriptstyle +$}}
  \put(242,135){\mbox{$\scriptstyle +$}}
  \put(221,89){\mbox{$\scriptstyle +$}}
  \put(236,87){\mbox{$\scriptstyle -$}}
  \put(263,89){\mbox{$\scriptstyle +$}}
  \put(248,87){\mbox{$\scriptstyle -$}}
  \put(195,-1){\mbox{$a$}}
  \put(291,-3){\mbox{$b$}}
  \put(315,-1){\mbox{$x$}}
        \end{picture}}
$$
Then applying the Proposition 4.2 from \cite{CP} we get the following
representative braid
$$
  \mbox{\begin{picture}(200,80)(0,0)
  \put(0,0){\epsfysize=80pt \epsfbox{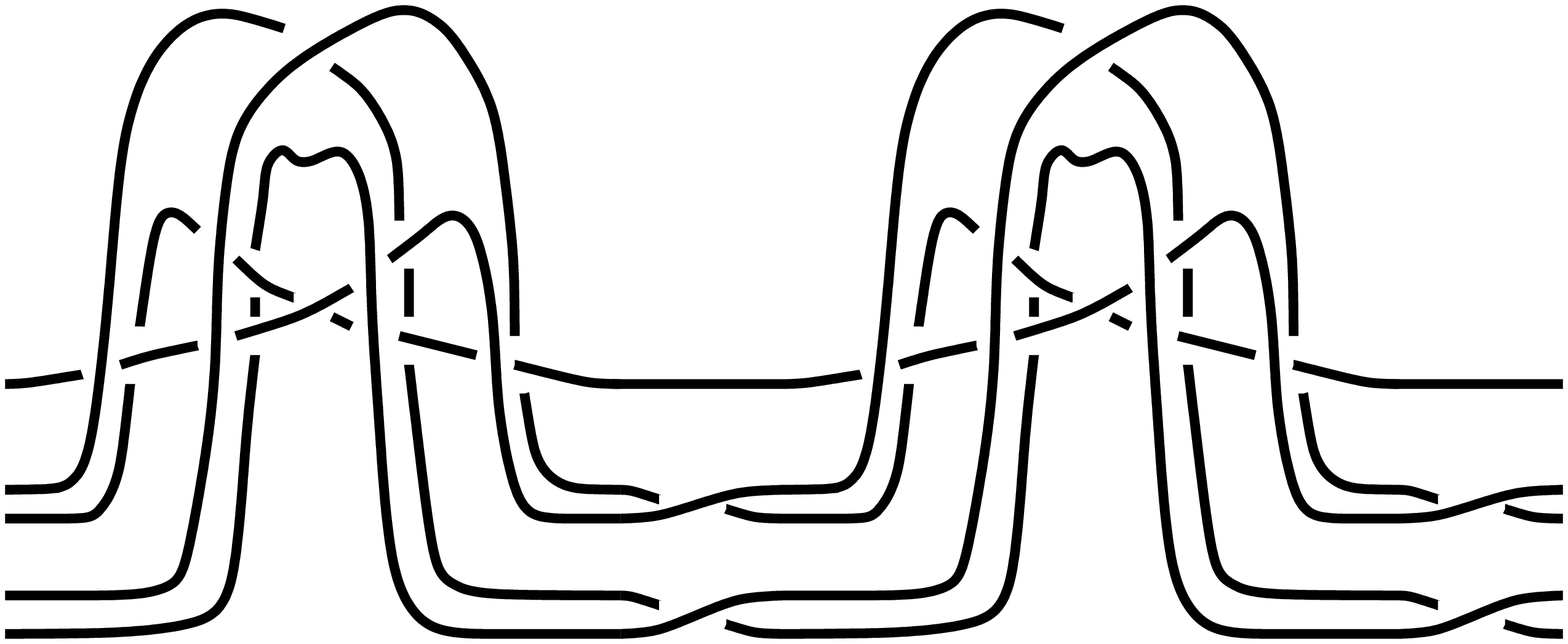}}
        \end{picture}}
$$
The closure of this braid gives a knot diagram isotopic to ours. It is
clear that the same arguments work for any generic divide as well.

\section*{Acknowledgments} I would like to express my deep gratitude
to the Department of Mathematics at Mount Holyoke College where
this work was done, for their warm hospitality and stimulating
atmosphere. I thank N.A'Campo, M.Hirasawa for valuable
comments.

\bibliographystyle{amsplain}

\vskip 1truepc
\noindent
\it
The Ohio State University,\\
1680 University Drive,\\
Mansfield, OH 44906\\
{\rm E-mail:} {\tt chmutov@math.ohio-state.edu}

\end{document}